\documentclass[11pt]{amsart}
\usepackage{amsmath, amsthm, amscd, amssymb}

\newtheorem{thm}{Theorem}[section]
\newtheorem*{thm*}{Theorem}
\newtheorem{cor}[thm]{Corollary}
\newtheorem*{cor*}{Corollary}
\newtheorem{lem}[thm]{Lemma}
\newtheorem{prop}[thm]{Proposition}

\newtheorem*{con*}{Conjecture}
\newtheorem*{prob*}{Problem}

\theoremstyle{definition}
\newtheorem{defn}[thm]{Definition}

\theoremstyle{remark}
\newtheorem{rem}[thm]{Remark}

\newcommand{\C}{\mathcal{C}}

\newcommand{\cL}{\mathcal{L}}

\newcommand{\bbC}{\mathbb{C}}

\newcommand{\bbR}{\mathbb{R}}

\newcommand{\val}{{\rm val}}
\newcommand{\im}{{\rm Im}}

\newcommand{\zigzag}{\mathbin{\raisebox{.2ex}{
      \hspace{-.4em}$\bigcirc$\hspace{-.75em}{\rm z}\hspace{.15em}}}}

\begin{document}
\title{Properties of the Generalized Zig-Zag Product of Graphs}

\author[Samuel Cooper]{Samuel Cooper$^{*}$}
\address{Department of Mathematics,
Vanderbilt University, Nashville, TN 37240}
\email{samuel.d.cooper@vanderbilt.edu}

\author[Dominic Dotterrer]{Dominic Dotterrer$^{*}$}
\address{Department of Mathematics, University of California, Santa Cruz, CA 95064} \email{ddotterr@ucsc.edu}
\thanks{$^{*}$Partially supported by an NSF REU grant}

\author[Stratos Prassidis]{Stratos Prassidis$^{**}$}
\address{Department of Mathematics
Canisius College, Buffalo, NY 14208, U.S.A.}
\email{prasside@canisius.edu}
\thanks{$^{**}$Partially supported by Canisius College Summer
Grant and an NSF REU grant}
%\date{07--04--2006}

\begin{abstract}
The operation of zig-zag products of graphs is the analogue of the semidirect product of groups. Using this observation, we present a categorical description of zig-zag products in order to generalize the construction for the category of simple graphs.  Also, we examine the covering properties of zig-zag products and we utilize these results to estimate their spectral invariants in general.  In addition, we provide specific spectral analysis for some such products.
\end{abstract}

\maketitle

\section{Introduction}

Expander graphs have very strong connectivity properties. For this reason,
they play a very important role in computer science, network design, 
cryptography,  just to mention a few of their applications. Originally, techniques from number theory provided the only known methods for constructing expanders. (\cite{lps}).
Zig-zag products of graphs were introduced as a combinatorial way of constructing infinite
sequences of expanders (\cite{alw}, \cite{rvw1}, \cite{rvw2}).  We begin by extending zig-zag products to the whole category of simple graphs (not necessarily finite).  We continue by showing that the spectral analysis of these graphs is greatly simplified by the structure of the zig-zag product for both finite and infinite graphs.

Zig-zag products were introduced as the graph-theoretic analogues of semidirect products 
of groups. That means that, with the right choice of generating sets,
the Cayley graph of a semidirect product of groups is the zig-zag product
of the Cayley graphs of those groups. In the first part of the paper we exploit 
this connection to give a categorical description of zig-zag products. That 
allows us to generalize the zig-zag construction for any graphs, finite
or infinite, regular or not. In general, we construct the zig-zag 
product $G{\zigzag}H$ when $G$ is $H$-{\it labeled} i.e., equipped
with a function
$${\alpha}: D(G) \to V(H)$$
where $(u, e)\in D(G)$ means that $e$ is an edge in $G$ with one end-point
$u$.

We begin by considering {\it locally constant} labelings and derive some specific information about the spectrum of the resulting zig-zag product.  Namely, the non-trivial eigenfunctions of the adjacency operator correspond bijectively between $G$ and $G \zigzag H$ with eigenvalues scaled only by a fixed constant.

Kesten's Theorem (\cite{ke1}, \cite{ke2}) implies that, for a $k$-regular
graph
$$\frac{2\sqrt{k-1}}{k} \le {\rho}(G) \le 1$$
with ${\rho}(G)$ maximum if and only if $G$ is amenable and minimum if and only if $G$ is the
$k$-regular tree.
In particular, using Kesten's Theorem, that if $G$ is amenable so is 
$G{\zigzag}H$. In this case, we construct a F{\o}lner sequence
for $G{\zigzag}H$ from a F{\o}lner sequence of $G$.

Covers can be used to estimate the spectra of graphs.  First we show that graph covers and combinatorial covers are preserved under generalized zig-zag products. Also, we show that for locally constant $H$-labelings of $G$ the zig-zag product $G{\zigzag}H$ is itself a combinatorial cover of $G$.  Combining these results, we construct an infinite sequence of combinatorial covers for any simple graph $G$.  

More specifically, 
let $G,H$ be locally finite simple graphs.  Let $\alpha$ be a locally constant $H$-labeling of G such that for all $h \in \im(\alpha), \; \val(h)=m$.
Define
\begin{enumerate}
\item $\Delta_1 = G \text{  and  } \alpha_1 = \alpha,$
\item $\Delta_{n+1} = \Delta_n \zigzag_{\alpha_n} H \text{ and } \alpha_{n+1}=\alpha_n \circ \tilde{\pi}_{n+1}.$
\end{enumerate}
Where $\tilde{ \pi }_{n+1}:D(\Delta_{n+1}) \rightarrow D(\Delta_n)$ is the map induced by the graph epimorphism,
$$\pi_{n+1} : \Delta_{n+1} \rightarrow \Delta_{n}.$$
Then $\pi_{n+1}$ is a combinatorial covering map for $n \geq 0.$ 
In fact, the spectrum of $\Delta_n$ is just the spectrum of $\Delta_{n-1}$ scaled by a constant $m$, fixed for all $n$.  

By applying the results of the previous section, we extend the result that eigenfunctions with non-zero eigenvalues correspond bijectively between $\Delta_n$ and $\Delta_{n-1}$.
As a consequence of this, given a graph exhibiting a spectral gap, it is possible to construct an infinite sequence of graphs of increasing complexity that also exhibit a spectral gap and, in fact, the length of the spectral gap is scaled by a fixed natural number for each term in the sequence, allowing for an arbitrarily large spectral gap.

The first two authors would like to thank Terry Bisson for his helpful
and insightfull suggestions. They would like also to thank  Canisius College
for the hospitality during the R.E.U. program in Summer 2006.

\section{Definitions and Notation}

All graphs will be simple locally finite graphs.
Let $G$ be a {\it simple graph} i.e., a graph without multiple edges and loops. 
Let $V({G})$ denote the vertex set of $G$ and
$E({G})$ the edge set. 
For two adjacent vertices $u$ and $v$ in $G$, 
(denoted $u \sim_G v$) we write $\{ u,v \}$ for the edge between $u$ and $v$, and $E_{G}(u)$ for the set
of edges adjacent to $u$. Also, we denote the neighborhood of a vertex $u \in V(G)$ by
$N_G(u)= \{ v \in V(G) : v \sim_G u \}$.
Set
$$D(G) = \{(u, e): u \;\text{is an end point of}\; e\}\subset
V(G) {\times} E(G).$$

Given two graphs, $G$ and $G'$, we call the map $\phi : V(G) \to V(G')$ a $\it{graph \; morphism}$ if,
$$u \sim_G v \Rightarrow \phi(u) \sim_{G'} \phi(v)$$
A graph morphism $\phi$ as above induces a map:
$$D({\phi}): D({G}) \to D({G}'), \; D({\phi})(u, \{u, v\}) =
({\phi}(u), \{{\phi}(u), {\phi}(v)\}).$$

Let $f, g: A \to {G}$ be two maps from a set $A$ to a graph $G$. Then
$f$ and $g$ are called {\it adjacent} if $f(a) \sim_{G} g(a)$ for
all $a\in A$.

Let $H$ be a graph. An  $H$-labeling on a graph $G$ is 
a function
$${\alpha}: D(G) \to V(H).$$
The pair $(G, {\alpha})$ is called an $H$-labeled graph.

\begin{defn}
Let $(G, {\alpha})$ be an $H$-labeled graph, with $H$ non-trivial. The zig-zag product
$G{\zigzag}H$ is the graph defined as follows:
\begin{itemize}
\item $V(G{\zigzag}H) = \{ (u,i) \mid \text{ there exists }\{ v,u \} \in E(G) \text{ such that } i \in N( \alpha (u, \{ v,u \} )) \}$.
\item $\{ (u,i),(v,j) \} \in E(G \zigzag H)$ if there is an edge $e= \{ u,v \} \in E(G)$
such that $i \sim {\alpha}(u,e)$ and $j \sim {\alpha}(v,e)$.
\item we denote the edge $\mathbf{e}=(e, \epsilon_1, \epsilon_2)$ where $\epsilon_1= \{ i, \alpha (u,e) \}$ and $ \epsilon_2 = \{ j, \alpha (v,e) \}$
\end{itemize}
\end{defn}

\begin{rem} 
The definition given generalizes the classical definition in several different
directions.
\begin{enumerate}
\item We do not assume that $G$ and $H$ are finite graph as in the classical
definitions of zig-zag products (\cite{alw}, \cite{rvw1}, \cite{rvw2}).
It is clear that if $G$ and $H$ are locally finite, then so is $G{\zigzag}H$. 
\item We do not require the graphs involved to be regular. The 
disadvantage is that $G{\zigzag}H$ is not necessarily regular.
\item The definition given above also generalizes the classical definitions 
in a different direction. In the
zig-zag products given in \cite{rvw1} and \cite{rvw2} the maps 
${\alpha}(u, -)$
are assumed to be bijections, for each $u \in V(G)$ and in \cite{alw} the maps
${\alpha}(u, -)$ are assumed to be injections whose image has $k$ elements,
where $k \mid {|V(G)|}$.
\end{enumerate}
\end{rem} 

\begin{rem}
\begin{enumerate}
\item In the usual definition of the zig-zag product, $V(G{\zigzag}H)$ was taken to be $V(G){\times}V(H)$. Our definition modifies this convention to avoid the presence of isolated vertices in $G{\zigzag}H$. 
\item There always exists a graph epimorphism,
$$\pi :G \zigzag H \rightarrow G, \quad (u,i) \mapsto u$$
\end{enumerate}
\end{rem}

We begin with some helpful preliminaries regarding the structure of the zig-zag product.

\begin{lem}\label{val} Let $G$ and $H$ be locally finite simple graphs with $\alpha$ an H-labeling of G.
Then for $(u,i) \in G {\zigzag}_{ \alpha }H,$

$${\val}(u,i)= \sum_{\{v \in V(G) \, \mid \, \alpha (u, \{ v,u \} ) \sim i\}} {\val}( \alpha (v, \{ v,u \} )).\ $$
\end{lem}

\begin{proof}
Notice that
$(v,j) \sim (u,i)$ in $G \zigzag_{\alpha} H$ if and only if there is an edge $\{v,u\} \in E(G)$ such that $\alpha (v, \{ v,u \} ) \sim j$ and $\alpha (u, \{ v,u \} ) \sim i.$  Now, for a fixed $v \in V(G)$ the number of such $j$ is equal to the ${\val}( \alpha (v, \{ v,u \} )$.  We count this value for each such $v \in V(G)$ such that $\alpha (u, \{ v,u \} ) \sim i $.
\end{proof}

\begin{lem}
Let $G$ be a finite graph,  $H$ a locally finite graph, and $\alpha$ an $H$-labeling of $G$.  Then
$$|E(G {\zigzag}_{ \alpha } H)|= \sum_{e =\{u, v\} \in E(G)}{\val}( \alpha (u,e)) \cdot {\val}( \alpha (v,e)).$$ 
\end{lem}

\begin{proof}
For $e = \{u, v\} \in E(G)$, let ${\val}(\alpha (u,e))=m,$ and ${\val} ( \alpha (v,e))=n$.  Since $H$ is simple, there exist distinct $i_1 , ... , i_l \in V(H) \ $ and distinct $j_1 , ... , j_m \in V(H)$ such that $i_k \sim \alpha (v,e), j_l \sim \alpha (u,e).$  Thus $ (u, j_l) \sim (v, i_k ) \ $ in $G{\zigzag}_{\alpha}H$, for $1 \leq k \leq n,$ and $1 \leq l \leq m$, which implies that $e$ induces $mn$ edges in $G {\zigzag}_{ \alpha } H$
\end{proof}

We now define a labeling that we will find particularly useful.  

\begin{defn}
Let G and H be locally finite simple graphs.  Then $\alpha: D(G) \rightarrow V(H)$ is called
locally constant if for each $u \in V(G)$, $\alpha(u, \{v, u\}) = \alpha(u, \{w,u\})$, for all $v$ and $w$ adjacent to $u$.
\end{defn}

\begin{lem}
Let $G$ be a finite graph, $H$ be locally finite graph and let $\alpha$ be a locally constant $H$-labeling such that for all $h \in {\rm Im}(\alpha)$, ${\val}(h)=k$.  Then $G \zigzag H$ contains at least $k^{|V(G)|}$ subgraphs isomorphic to $G$.
\end{lem}

\begin{proof}
For each $u \in V(G)$ choose $(u,i) \in V( G \zigzag H )$.  There are $k$ such choices for each $u$.  Let $G'$ be the subgraph spanned by the chosen vertices.  Note that $(u,i) \sim (v,j) \Leftrightarrow v \sim_G u$.  Thus, since a bijection exists between $G'$ and $G$, they are isomorphic.  Note that there exists $k^{|V(G)|}$ choices for the construction of $G'$.
\end{proof}

\section{Categories of $H$-labeled Graphs} 

We will give a categorical description of the zig-zag 
product. We define two categories of $H$-labeled graphs. 
In one the morphism are commutative diagrams ``up to adjacency'' 
and in the other the diagrams strictly commute.

We write ${\cL}_w(H)$ for the {\it weak category} of $H$-labeled graphs:
\begin{itemize}
\item Objects are $H$-labeled graphs $(G, {\alpha})$.
\item A morphism ${\phi}: (G_1, {\alpha}_1) \to (G_2, {\alpha}_2)$ is
a graph map ${\phi} = ({\phi}, \{{\phi}_{u,v}\})$ such that:
$${\alpha}_1(u, \{u,v\}) \sim h \;\Longrightarrow\; 
{\alpha}_2({\phi}(u), {\phi}_{u,v}(\{u,v\})) \sim h, 
\;\;\text{for}\;\; h \in V(H),$$
where ${\phi}_{u,v}(\{u,v\}) =\{{\phi}(u), {\phi}(v)\}.$
\end{itemize}
We also write ${\cL}_s(H)$ for the {\it strict category} of 
$H$-labeled graphs, with the same objects but morphisms maps
$${\phi}: 
(G_1, {\alpha}_1) \to (G_2, {\alpha}_2)$$
such that the following diagram commutes:
$$\begin{CD}
D(G_1) @>{{\alpha}_1}>> V(H) \\
@V{D({\phi})}VV @VV{\text{id}}V \\
D(G_2) @>{{\alpha}_2}>> V(H)
\end{CD}$$ 

There is a forgetful functor:
$${\iota}: {\cL}_s(H) \to {\cL}_w(H).$$

We now describe the dependence of the zig-zag product on the second variable.

\begin{lem}\label{lem-induce}
Let ${\psi}: H_1 \to H_2$ be a graph morphism.
\begin{enumerate}
\item The morphism $\psi$ induces a functor:
$${\Psi}_s: {\cL}_s(H_1) \to {\cL}_s(H_2).$$
\item Assume that ${\psi}$ is onto on the set of vertices and satisfies:
$$h \sim_{H_1} h' \;\Longleftrightarrow\; {\psi}(h) \sim_{H_2} {\psi}(h').$$
Then $\psi$ induces a functor 
$${\Psi}_w: {\cL}_w(H_1) \to {\cL}_w(H_2).$$
\end{enumerate}
\end{lem}

\begin{proof}
If $(G, {\alpha})$ is an $H_1$-labeled graph, then 
$(G, {\psi}{\circ}{\alpha})$ is an $H_2$-labeled graph. Part (1) follows 
immediately. For Part (2) and for a morphism
$$({\phi}, \{{\phi}_{u,v}\}): (G, {\alpha}_1) \to (G, {\alpha}_2),$$ 
in ${\cL}_w(H_1)$, notice that:
$$\begin{array}{rlll}
{\psi}({\alpha}_1(u, \{u, v\})) &  \sim_{H_2}  h & \Longrightarrow & \\
{\psi}({\alpha}_1(u, \{u, v\})) &  \sim_{H_2}  {\psi}(h') 
& \Longrightarrow & (\psi \;\;\text{is onto})\\
{\alpha}_1(u, \{u, v\}) &  \sim_{H_1} h' 
& \Longrightarrow & (\text{by the property of $\psi$}) \\
{\alpha}_2({\phi}(u), {\phi}_{u,v}(\{u, v\})) &  \sim_{H_1}  h' 
& \Longrightarrow & (\text{by the property of $\phi$}) \\
{\psi}({\alpha}_2({\phi}(u), {\phi}_{u,v}(\{u, v\}))) &  
\sim_{H_2}  {\psi}(h') 
& \Longrightarrow & (\text{since $\psi$ is a graph morphism} ) \\
{\psi}({\alpha}_2({\phi}(u), {\phi}_{u,v}(\{u, v\}))) &  \sim_{H_2}  h && 
\end{array}$$
proving that 
$$({\phi}, \{{\phi}_{u,v}\}): (G, {\psi}{\circ}{\alpha}_1) \to 
(G, {\psi}{\circ}{\alpha}_2),$$ 
is a morphism in ${\cL}_w(H_2)$.
\end{proof}

We now describe the dependence of the zig-zag product on the second variable.
Again, we have two different descriptions, one for each category.
Let ${\C}_s$  be a category with
\begin{itemize}
\item Objects are triples $((G, {\alpha}), H)$ where $(G, {\alpha})$ is in
${\cL}_s(H)$, 
\item A morphism 
$$({\phi}, {\psi}): ((G_1, {\alpha}_1), H_1) \rightarrow ((G_2, {\alpha}_2), H_2)$$
is a pair of morphisms:
$$({\phi}, {\psi}) \in \text{Mor}_{{\cL}_s(H_1)}((G_1, {\alpha}_1),
(G_2, {\alpha}_2)){\times} \text{Mor}(H_1, H_2).$$
\end{itemize}

Similarly, we define ${\C}_w$ to be the category with the same set of
objects as ${\C}_s$ and morphisms pairs $({\phi}, {\psi})$ as above where
$\phi$ is a morphism in ${\cL}_w(H_1)$ and $\psi$ satisfies Condition (2)
in Lemma \ref{lem-induce}.

If a pair $({\phi}, {\psi})$ is a morphism in ${\C}_s$ or ${\C}_w$. Define
$$F : V(G_1{\zigzag}H_1) \to V(G_2{\zigzag}H_2), \;\;
F(u, h) = ({\phi}(u), {\psi}(h)).$$
Also, F induces a map:
$$E(G_1 \zigzag_{\alpha_1} H) \rightarrow E(G_2 \zigzag_{\alpha_2} H), \;\; \{(u,i),(v,j)\} \mapsto (\{ \phi(u),\phi(v) \}, \{ \psi(\alpha_1 (u,\{u,v\}), \psi(i)\},  \{ \psi(\alpha_1 (v,\{u,v\}), \psi(j) \} )$$

\begin{prop}\label{prop-morphism}
The pair induces a graph map:
$$F : V(G_1{\zigzag}H_1) \to V(G_2{\zigzag}H_2).$$
\end{prop}

\begin{proof}
It follows from the definitions.
\end{proof}

Thus combining the categorical properties of the construction,
we will find conditions that determine when zig-zag products are isomorphic.

\begin{thm}\label{thm-isomorphism}
Let ${\phi}: G_1 \to G_2$ and ${\psi}: H_1 \to H_2$ be graph isomorphisms.
Also, assume that $(G_i, {\alpha}_i)$ is a $H_i$-labeled graph ($i = 1, 2$)
\begin{enumerate}
\item If $({\phi}, {\psi})$ is a morphism in ${\C}_s$, then 
$$F: G_1{\zigzag}H_1 \to G_2{\zigzag}H_2$$
is a graph isomorphism.
\item If $({\phi}, {\psi})$ is a morphism in ${\C}_h$ that satisfies
$${\alpha}_1(u, \{u, v\}) \sim_{H_1} h \;\Longleftrightarrow\;
{\alpha}_2({\phi}(u), {\phi}_{u,v}(\{u, v\}) \sim_{H_2} h.$$
 Then
$$F: G_1{\zigzag}H_1 \to G_2{\zigzag}H_2$$
is a graph isomorphism.  
\end{enumerate}
\end{thm}

\begin{proof}
The conditions imposed guarantee that there is a natural inverse graph map 
to $F$.
\end{proof}

The next result is another graph theoretic analogue of the semidirect
product construction on groups. Let $(G, {\alpha})$ be an $H$-labeled
graph, and $Z$ a graph equipped with an graph map $f: Z \to G$. Then
$f$ induces a map
$${\alpha}{\circ}D(f): D(Z) \to D(G) \to V(H).$$
Also, if $g: V(Z) \to V(H)$ is a  map, where $H$ is a graph, define 
$$g{\circ}\text{pr}:D(Z) \to V(Z) \to V(H).$$
The following follows from the definitions.

\begin{prop}\label{prop-product}
Let $Z$ be a graph, $f: Z \to G$ a graph map and $g: V(Z) \to V(H)$ a map 
such that
${\alpha}{\circ}D(f)$ is adjacent to $g{\circ}\text{pr}$. Then there exists an edge ${\epsilon}_{z,e}$ in H, that joins ${\alpha}{\circ}D(f)(z,e)$ to
$g{\circ}\text{pr}(z, e)$.
Then the pair 
$(f, g)$
$$(f, g): V(Z) \to V(G{\zigzag}H), \;\; (f, g)(z) = (f(z), g(z)),$$
and, for $z \sim_Z w$,
$$(f,g)(\{z,w\})= \{(f(z),g(z)),(f(w),g(w))\}$$
\end{prop}

\section{Spectral Analysis of  Zig-Zag Products}

For a graph $G$, the {\it adjacency matrix} of $G$ is the matrix
$A$ with rows and columns indexed by the vertex set of $G$ such that
$$A(u,v) = \left\{
\begin{array}{cl}
1, & \text{if}\; u \sim v \\[2ex]
0, & \text{otherwise}
\end{array}\right.$$

\begin{defn}
Let $L^2(G)= \{ f: V(G) \rightarrow \bbC \mid \sum_{v \in V(G)} |f(v)|^2 < \infty \}$
Define the adjacency operator on $L^2(G)$:
$$
A:L^2(G) \rightarrow L^2(G), \quad
f \mapsto Af, \quad Af(v) = \sum_{u\sim v}f(u).$$
\end{defn}

By the {\it eigenvalues of} $G$ we mean the eigenvalues of $A$.
We write
$${\rho}(G) = \max\{|{\lambda}|: \;\text{an eigenvalue of}\; A \}$$
called the  {\it spectral radius} of $G$, and
$${\lambda}(G) = \max\{|{\lambda}|: \;\text{an eigenvalue of}\; A \;
|{\lambda}| \neq \rho (G) \}.$$

\begin{thm}
Let $G$ be an $H$-labeled graph. Assume that the labeling $\alpha$ is locally constant with the 
valency of $h$ equal to $n$ for all $h$ in the image of $\alpha$.
 Let $f \in L^2(G)$ be an eigenvector of $A_G$ with eigenvalue 
 $ \lambda \ $.
Then there exists an eigenvector of $ A_{G \zigzag H} $ with eigenvalue $k \lambda$.
\end{thm} 

\begin{proof}
Let $\hat{f} \in L^2(G {\zigzag} H)$ be defined as $ \hat{f} (u,i)=f(u)$, for all 
$(u,i) \in V(G {\zigzag} H)$. Notice that for each $u \in V(G)$,
$$\sum_{i \in V(H), (u, i) \in  V(G{\zigzag}H)}\hat{f}(u, i) = nf(u).$$ 
We claim that $\hat{f}$ is the desired eigenvector.
First we show that $ \hat{f} \in L^2(G \zigzag H)$.  This is true because 
$$ \sum_{(u,i) \in V(G {\zigzag}_{ \alpha } H)} | \hat{f} (u, i) |^2 =n \sum_{u \in V(G)} |f(u)|^2 < \infty.  $$
In particular,
$$\| \hat{f}\|_{L^2(G{\zigzag}H)} = \sqrt{n}\| f\|_{L^2(G)}.$$
Now, since $f$ is an eigenvector of $A_G$, we know that
$$ \lambda f(u)= A_Gf(u) = \sum_{v \sim u} f(v)$$ for all $u \in V(G).$
But this implies by Lemma \ref{val} that 
$$A_{G \zigzag H} \hat{f} (u,i)  =  \sum_{(v,j) \sim (u,i)} \hat{f} (v,j)
= \sum_{v \sim u} \quad \sum_{j \sim \alpha(v, \{ v,u \} ) } \hat{f} (v,j)  
= n \sum_{ v \sim u} f(v),$$
since $\alpha$ being locally constant implies that $ \{ v \mid v \sim u \} = \{ v \mid v \sim u, \alpha (u, \{ v,u \} ) \sim i \} \ $.  Thus
$$ A_{G \zigzag H} \hat{f} (u,i) = n \sum_{ v \sim u} f(v) = k \lambda f(u)
= n \lambda \hat{f} (u,i) \ ,$$
 for all $(u,i) \in V(G {\zigzag}_{ \alpha } H))$.
\end{proof}

\begin{thm}
Let $G,H,\alpha$ be as above.
Let $ \hat{f} $ be an eigenvector of $ A_{G \zigzag_{\alpha} H}$ with eigenvalue $
\lambda$.
Then either $\lambda =0$ or $ { \lambda}/n $ is an eigenvalue of
$ A_{G}$.
\end{thm}

\begin{proof}
From the definition of the zig-zag product we can see that $$\sum_{(v,j)
\sim (u,i)} {\hat{f} (v,j)} = \sum_{\substack{v\sim u \\ \alpha (v, \{u,
v\}) \sim j \\ \alpha (u, \{u, v\}) \sim i}} {\hat{f}(v, j)} = \sum_{v \sim u}
\sum_{\alpha(u, \{u, v\}) \sim i} {\hat{f} (v,j)}.$$  
the last equality holds since
$\alpha$ is locally constant. Thus, $$\sum_{(v,j)
\sim (u,i)} {\hat{f}(v,j)} = \sum_{(v,j) \sim (u,k)} {\hat{f}(v,j)},$$
and since $\hat{f}$ is an eigenvector, $$\lambda \hat{f}(u,i) = \lambda \hat{f} (u,k).$$ 
Thus, if $\lambda$ is not equal to $0$, then $$\hat{f} (u,i) = \hat{f}
(u,k).$$
Thus we see that $\hat{f}(u,i)$ depends only on $u$.  Now, define $f \in
L^2(G)$ by the formula $$f(u) = \hat{f} (u,i),$$ for all $i \in V(H)$. 
So, $$Af(v) = \sum_{u \sim v} {f(u)} = \sum_{u \sim v} {\hat{f}(u,i)}.$$ 
Now, since for any $u, v \in V(G)$ such that $u \sim v$, there are $n$
vertices $i \in V(H)$ such that $$i \sim \alpha(u, \{u,v\}),$$
we see that $$\sum_{u \sim v} {f(u,i)} = \frac{1}{n} \sum_{\substack{{u
\sim v} \\ i \sim {\alpha(u, \{u,v\})}}}{\hat{f}(u,i)}.$$
So, for any $j \in V(H)$, this last sum is equal to $$\frac{1}{n} (\lambda
\hat{f}(v,j)) = \frac{\lambda}{n}f(v).$$
\end{proof}

The normalized adjacency matrix $P = (1/d)A_G$ of a $d$-regular graph $G$ can be 
considered as an operator on $L^2(V(G))$. Then the spectral radius of $G$
is given by
$${\rho}_N(G) = \| P \| = \sup_{f\in L^2(V(G)) \setminus \{0\}}\frac{
{\langle}Pf, f{\rangle}}{\| f \|}.$$
For a $d$-regular graph,
$$Pf(v) = \frac{1}{d}\sum_{u \sim v}f(u).$$

\begin{thm}\label{thm-radius}
Let $G$ be a $m$-regular graph, $H$ a $d$-regular graph on $m$ vertices.
Let $\alpha$ be an $H$-labeling such that ${\alpha}(u, -)$ is a
bijection for each $u \in V(G)$. Then
$${\rho}_N(G) \le {\rho}_N(G{\zigzag}H).$$
\end{thm}

\begin{proof}
Let $f\in L^2(G)$. Then $f$ induces an element $\hat{f}\in L^2(G{\zigzag}H)$ by
$$\hat{f}: G{\zigzag}H \to \bbC, \quad \hat{f}(u, i) = f(u).$$
It is obvious that $\hat{f}\in L^2(G{\zigzag}H)$. Then
$$\| \hat{f} \|_z = \sum_{(v,i)\in V(G{\zigzag}H)}|\hat{f}(v,i)|^2 = 
k\sum_{v\in V(G)}|f(v)|^2 = \| f\|,$$
where $\| - \|_z$ denotes the norm in $L^2(G{\zigzag}H)$.

Let
$v\in V(G)$ and  $u\in N(v)$. Let ${\alpha}(u, |\{u, v\}) = i$ and
${\alpha}(v, |\{u, v\}) = j$. Then $(v, m) \sim (u, n)$ in
$G{\zigzag}H$ if and only if $m \sim i$ and $n \sim j$. Since $H$ is
$d$-regular, there are $d^2$ such pairs. If $Q$ denotes the normalized
adjacency matrix in $G{\zigzag}H$, then, if ${\langle}-, -{\rangle}_z$ denote
the inner product in $L^2(G{\zigzag}H)$,
$$\begin{array}{lll}
{\langle}Qf, f{\rangle}_z & = & \displaystyle{\frac{1}{d^2} 
\sum_{(v,i)\in V(G{\zigzag}H)}\left(\sum_{(u,j)\sim (v,i)}\hat{f}(u,j)\right)
\overline{f'(v, i)}} \\[3ex]
& = & \displaystyle{
\frac{1}{d^2} \sum_{(v,i)\in V(G{\zigzag}H)}
\left(\sum_{(u,j)\sim (v,i)}f(u)\right)\overline{f(v)}} \\[3ex]
& = & \displaystyle{
\frac{1}{d^2} \sum_{v\in V(G)}
\left(\sum_{\stackrel{(u,j)\sim (v,i)}{\text{for some $i$, $j$}}}
f(u)\right)\overline{f(v)}} \\[3ex]
& = & \displaystyle{\frac{1}{d^2}\sum_{v\in V(G)}
\left(\sum_{u\sim v}\sum_{\stackrel{(u,j)\sim (v,i)}{\text{for some $i$, $j$}}}
f(u)\right)\overline{f(v)}} \\[3ex]
& = & \displaystyle{\frac{1}{d^2}\sum_{v\in V(G)}
\left(\sum_{u\sim v}d^2f(u)\right)\overline{f(v)}} \\[3ex]
& = & \displaystyle{\sum_{v\in V(G)}
\left(\sum_{u\sim v}f(u)\right)\overline{f(v)}} \\[3ex]
& = & k{\langle}Pf, f{\rangle}
\end{array}
$$
Using this calculation, we get
$$\begin{array}{lll}
{\rho}(G{\zigzag}H) & = &\displaystyle{
\sup_{g\in L^2(G{\zigzag}H) \setminus \{0\}}\frac{
{\langle}Qg, g{\rangle}_z}{\| g \|_z} \ge
\sup_{f\in L^2(G) \setminus \{0\}}\frac{
{\langle}Q\hat{f}, \hat{f}{\rangle}_z}{\| \hat{f}\|_z}} \\[2ex]
& = & \displaystyle{\sup_{f\in L^2(G) \setminus \{0\}}\frac{
k{\langle}Pf, f{\rangle}}{k\| f \|}
= \sup_{f\in L^2(G) \setminus \{0\}}\frac{
{\langle}Pf, f{\rangle}}{\| f \|} = {\rho}(G)}
\end{array}
$$
\end{proof}

Theorem \ref{thm-radius} and Kesten's Theorem imply that the zig-zag product
$G{\zigzag}H$ of graphs that satisfy the assumptions of Theorem 
\ref{thm-radius} with $G$ amenable, is amenable.

If $G$ has enough symmetries then
(\cite{day}, \cite{ke1}, \cite{ke2}, \cite{woess})
$$\frac{2\sqrt{d - 1}}{d} \le {\rho}(G) \le 1$$
and equality on the left was obtained if and only if $G$ is a $d$-regular
tree and on the right if and only if $G$ is amenable. We recall the definition
of amenability that we will use:

\begin{defn}
A graph $G$ is amenable if there is a sequence of finite subgraphs 
$\{F_n\}_{n=1}^{\infty}$ whose union is $G$ such that
$$\lim_{n\to \infty}\frac{|{\partial}F_n|}{|F_n|} =0.$$
Here, for a subgraph $H$ of $G$,
$${\partial}H = \{e\in E(G): \;\text{one endpoint of $e$ is in $H$ and the
other it is not}\}.$$
The sequence $\{F_n\}_{n=1}^{\infty}$ is called a F{\o}lner Sequence.
\end{defn}

For an $H$-labeled graph $(G, {\alpha})$, let  $F$ be a 
subgraph of $G$. Then
the restriction of $\alpha$ induces an $H$-labeling on $F$. Using this
observation we can construct F{\o}lner Sequences in the zig-zag products.

\begin{thm}\label{thm-isoperimetric}
Let $(G, {\alpha})$ be an $H$-labeled graph satisfying the hypotheses of Theorem \ref{thm-radius},
Let $\{F_n\}_{n=1}^{\infty}$ be a F{\o}lner Sequence in $G$. 
Then $\{F_n{\zigzag}H\}_{n=1}^{\infty}$ is a F{\o}lner Sequence in 
$G{\zigzag}H$.

In particular, if $G$ is amenable, then $G{\zigzag}H$ is amenable for each
$H$-labeling.
\end{thm}

\begin{proof}
First notice that for each $n$, $|F_n{\zigzag}H| = |F_n|{\times}|H|$.
Let $D$ be the maximum degree in $H$. 
The definition of the zig-zag product implies that:
$$|{\partial}(F_n{\zigzag}H)| \le D^2|{\partial}F_n|.$$
Thus, taking limits,
$$\lim_{n \to \infty}\frac{|{\partial}(F_n{\zigzag}H)|}{|F_n{\zigzag}H|} = 0.$$
\end{proof} 

\section{Covers of Zig-Zag Product}

Despite the heavily relaxed conditions in our generalized definition of the zig-zag product, we retain some of the combinatorial properties of the original graphs.

First we show that zig-zag products preserve graph coverings.

\begin{defn}
For two graphs, $G, \tilde{G}$, a map $p: \tilde{G} \rightarrow G$ is said to be a graph covering map if:
\begin{enumerate}
\item $x \sim_{ \tilde{G} } y \Rightarrow p(x) \sim_G p(y)$ i.e., $p$ is a graph morphism.
\item For each $x \in p^{-1}(u)$, $p \mid_{N(x)} :N(x) \rightarrow N(u)$ is a bijection.
\end{enumerate}
In this case $\tilde{G}$ is said to cover, or be a graph covering of $G$.
\end{defn}

\begin{rem}
Let $G$ be a group and $S$ a finite symmetric generating set of $G$. Let $H$ be a normal subgroup such that $S{\cap}H = {\emptyset}$ and the cosets of $S$ generate $G/H$ without repetitions. Then the natural projection induces a graph covering map between the corresponding Cayley graphs:
$$p: \text{Cay}(G, S) \to \text{Cay}(G/H, S).$$
\end{rem}

\begin{thm}\label{covering}
Let $H$ be a graph and $(G, {\alpha})$ be an $H$-labeled graph. 
Let $p: \tilde{G} \rightarrow G$ be a graph covering map. 
Then there exists an $H$-labeling, $ \beta$ on $\tilde{G}$,
such that, the natural map
$$\hat{p} : \tilde{G} \zigzag_{ \beta } H \rightarrow G \zigzag_{ \alpha } H, \;\; \hat{p}(x, i) = (p(x), i),$$
is a covering map.
\end{thm}

\begin{proof}
The covering $p$ induces a map, 
$$\tilde{p} : D( \tilde{G} ) \rightarrow D(G), \quad ( x, \{ x,y \} ) \mapsto (p(x), \{ p(x),p(y) \} ).$$
Since $\{ p(x),p(y) \} \in E(G)$, $\tilde{p}$ induces an $H$-labeling of $\tilde{G}$ by
${\beta} = \alpha \circ \tilde{p}$.  Then $p$ is a morphism in ${\cL}_s(H)$ and induces a graph map 
(Proposition \ref{prop-morphism})
$$ \hat{p} : \tilde{G} \zigzag_{ \beta} H \rightarrow G \zigzag_{ \alpha } H, \quad (x,i) \mapsto (p(x),i).$$
Now let $(v,j) \in N(\hat{p} (x,i))$. By definition
$${\alpha}(v, \{ v, p(x)\}) \sim_H j, \quad  {\alpha}(p(x), \{ v, p(x) \}) \sim_H i \qquad (*)$$
Also, since $v \sim_G p(x)$ and 
$p$ is a graph covering map, there exists a unique $y \in p^{-1}(v)$ such that $y \sim_{ \tilde{G} } x$.
Clearly $\hat{p} (y, j) = (v, j)$. But, using (*),
$${\beta}(x, \{x, y\}) = {\alpha}(p(x), \{p(x), p(y)\}) = {\alpha}(u, \{u, v\}) \sim_H i.$$
Similarly, by (*) again, ${\beta}(y, \{ x,  y\}) \sim_H j$. Thus $(x,i) \sim (y,j)$ in ${\tilde{G}}{\zigzag}_{\beta}H$. That shows that the map
$$\hat{p}|_{N(x, i)} : N(x, i) \rightarrow N(\hat{p}((x, i)$$
is onto. The uniqueness of the choice for $y$ shows that $\hat{p}|_{N(x, i)}$ is injective. So
$\hat{p}|_{N(x, i)}$ is a bijection.
\end{proof}

The following follows immediately from \cite{woess}, Lemma 11.4.

\begin{cor}
If $\tilde{G}$ is a graph covering of an $H$-labeled graph $G$, with the notation of Theorem \ref{covering}, 
$$ \rho(\tilde{G} \zigzag_{\beta} H) \leq \rho(G \zigzag_{\alpha} H)$$
\end{cor}

Chung and Yau (\cite{cy}) introduced a variant to graph coverings. We will show that zig-zag products preserve combinatorial coverings.

\begin{defn}
For two graphs, $G, \tilde{G}$, a map $\pi : \tilde{G} \rightarrow G$ is said to be a combinatorial covering map if:
\begin{enumerate}
\item There is $m \in {\bbR}^+ \cup \{\infty\}$ such that for every  $ \{ u,v \} \in E(G)$,
$$| \{ \{ x,y \} \in E( \tilde{G} ): \pi (x) = u, \pi (y) = v \}| = m$$
\item For  $x,y \in V( \tilde{G} )$ with $ \pi (x) = \pi (y)$, and $v \sim_G \pi (x)$, we have,
$$ |N(x) \cap \pi^{-1} (v)| = |N(y) \cap  \pi^{-1} (v)|.$$
\end{enumerate}
\end{defn}

\begin{defn}
For a finite graph $G$, define the normalized Laplacian,
$$
\cL(u,v) = \left\{
\begin{array}{cll}
1, & \text{if}\; u=v \\
-\frac{1}{\sqrt{\val(u) \val(v)}}, &\text{if}\; u \sim_G v \\
0, & \text{otherwise}
\end{array}\right.$$
If $G$ is $d$-regular, then $\cL = I - P$.
\end{defn}

\begin{thm}\label{combinatorial}
Let $G, \tilde{G} ,H$ be graphs with $\alpha $ be an $H$-labeling of $G$.
Let $p: \tilde{G} \rightarrow G$ be a combinatorial covering map.  Then there exists an $H$-labeling $\beta$ of $\tilde{G}$ such that the induced map 
$$\hat{p} : \tilde{G} \zigzag_{ \beta } H \rightarrow G \zigzag_{ \alpha } H,$$
is a combinatorial covering map,
\end{thm}

\begin{proof}
As in the proof of Theorem \ref{covering}, $p$ induces a labeling $\beta$ on $\tilde{G}$ and 
a graph morphism
$$ \hat{p} : \tilde{G} \zigzag_{ \beta} H \rightarrow G \zigzag_{ \alpha } H, \quad (x,i) \mapsto (p(x),i).$$

Let $ \{ (u,i), (v,j) \} \in E(G \zigzag_{\alpha} H)$.  Then,
$$\begin{array}{ll}
& \{ \{ (x,k),(x',k') \} \in E( \tilde{G} \zigzag_{ \beta} H) : \hat{p}(x,k)=(u,i), \;\; \hat{p}(x',k')=(v,j)\} = \\
&\{ \{ (x,i),(x',j) \} \in E( \tilde{G} \zigzag_{ \beta } H) : x \in p^{-1}(u), \;\; x' \in p^{-1}(v) \}.
\end{array}$$
Let $\{ (u,i),(v,j) \} \in E(G \zigzag_{ \alpha } H)$, $x \in p^{-1}(u)$, $x' \in p^{-1}(v)$. Then, as in the proof of Theorem \ref{covering},
$${\beta}(x, \{x, x'\}) \sim_H i \Longleftrightarrow {\alpha}(p(x), \{p(x), p(x')\}) \sim_H i 
\Longleftrightarrow {\alpha}(u, \{u, v\}) \sim_H i.$$
The last relation is satisfied because $(u, i) \sim_{G{\zigzag}_{\alpha}{H}} (v, j)$. Similarly, 
$${\beta}(x', \{x, x'\}) \sim_H j \Longleftrightarrow {\alpha}(v, \{u, v\}) \sim_H j.$$
Therefore,
$$(x, i) \sim_{\tilde{G}{\zigzag}_{\beta}H} (x', j) \Longleftrightarrow x \sim_{\tilde{G}}x'.$$
Thus,
$$\begin{array}{ll}
& | \{ \{ (x,k),(x',k') \} \in E( \tilde{G} \zigzag_{\beta } H) : \hat{p}(x,k) = (u,i), \;\; \hat{p}(x',k')=(v,j) \} | = \\
& | \{ \{ x,x' \} \in E(\tilde{G}) : p(x)=u, \; p(x')=v \} |=m 
\end{array}$$
where $m$ is the index of $p$.

Now let $(x,i) \in V(\tilde{G} \zigzag_{ \beta } H)$ with $\hat{p}(x, i)$ adjacent to $(u,j)$ in 
$G{\zigzag}_{\alpha}H$. Then $u \sim_G p(x)$ and
$${\alpha}(u, \{u, p(x)\}) \sim_H j, \;\;{\alpha}(p(x), \{u, p(x)\}) \sim_H i.$$

If $y \in N(x){\cap}p^{-1}(u)$ then $p(y) = u$ and $y \sim _{\tilde{G}} x$. Then
\begin{enumerate}
\item ${\beta}(y, \{x, y\}) \sim_H j \Longleftrightarrow {\alpha}(p(y), \{p(x), p(y)\}) \sim_H j
\Longleftrightarrow  {\alpha}(u, \{p(x), u\}) \sim_H j$, which holds from above.
\item ${\beta}(x, \{x, y\}) \sim_H i \Longrightarrow {\alpha}(p(x), \{p(x), p(y)\}) \sim_H i
\Longleftrightarrow  {\alpha}(p(x), \{p(x), u\}) \sim_H i$, which holds from above.
\end{enumerate}
Thus $(y, j)$ is adjacent to $(x, i)$ in $\tilde{G}{\zigzag}_{\beta}H$ and 
$(y, j) \in N(x, i)$. Also, $\hat{p}(y, j) = (p(y), j) = (u, j)$, which implies
$$(y, j) \in N(x, i){\cap}\hat{p}^{-1}(u, j).$$
 Thus there is a map
$${\xi}: N(x){\cap}p^{-1}(u) \to N(x, i){\cap}\hat{p}^{-1}(u, j), \quad y \mapsto (y, j).$$
But notice that 
$$N(x, i) \cap \hat{p}^{-1}(u, j) = \{ (z, k) : (z, k) \sim (x, i), \; \hat{p}(z, k)=(u, j) \} =
 \{ (z, j) : \; z\sim_{\tilde{G}}x, \; p(z) = u \}. $$
 So $\xi$ is a bijection and the two sets have the same cardinality.

Let $(x, i),(y, j) \in V(\tilde{G} \zigzag_{ \beta } H)$ such that $\hat{p}(x, i)=\hat{p}(y, j)$ and 
$ \hat{p}(x, i)$, $\hat{p}(y, j)$ are adjacent to $(u, j)$ in $G{\zigzag}_{\alpha}H$. That implies $i = j$. 
Using the above result
$$\begin{array}{llll}
|N(x, i) \cap \hat{p}^{-1}(u, j)| & = & |N(x){\cap}p^{-1}(u)|  &\\
&= &  |N(y){\cap}p^{-1}(u)| & \text{(because $p$ is a combinatorial cover)} \\
& = & |N(y, j) \cap \hat{p}^{-1}(u, j)| & \\
& = & |N(y, j) \cap \hat{p}^{-1}(u, j)| &
\end{array}$$
That completes the proof.
\end{proof}

In \cite{cy}, Lemma 1, the connection between the spectra of graphs in combinatorial covers is given. Applying this observation, we get the following consequence of Theorem \ref{combinatorial}.

\begin{cor}\label{laplace}
Let $\tilde{G}$, $G$ and $H$ be finite graphs.
Let $\tilde{G}$ be a combinatorial cover of an $H$-labeled graph $G$. With the notation of Theorem \ref{combinatorial}, each of the eigenvalues of the normalized Laplacian of $G{\zigzag}_{\alpha}H$ are eigenvalues of the normalized Laplacian of $\tilde{G}{\zigzag}_{\beta}H$. In particular, if both $G{\zigzag}_{\alpha}H$ and $\tilde{G}{\zigzag}_{\beta}H$ are regular,
$$\text{spec}(G{\zigzag}_{\alpha}H) \subset \text{spec}(\tilde{G}{\zigzag}_{\beta}H).$$
\end{cor}

Combinatorial graphs arise naturally in zig-zag products.

\begin{thm}
Let $G$ be a locally constant $H$-labeled graph. Assume that the $H$-labeling $\alpha$ is such that for all $h \in \text{Im}(\alpha), \; {\val}(h) = n$.  Then
$$\pi :G \zigzag H \rightarrow G, \quad (u,i) \mapsto u$$
is a combinatorial covering map.
\end{thm}

\begin{proof}
Let $\{ v,u \} \in E(G),$
Then $$| \{ \{ x,y \} \in E(G \zigzag H) \mid \pi (x)=u, \; \pi (y)=v\} | = | \{ \{ (v,j),(u,i) \} \in E(G \zigzag H) \} |=n^2.$$
Now consider $(u,i),(u,j) \in \pi^{-1} (u)$.  Let $v \sim u$
Then $$|N(u,i) \cap \pi^{-1} (v)| = | \{ (v,k) \mid k \sim \alpha (v, \{ u,v \} ) \} | = n =  |N(u,j) \cap \pi^{-1} (v)|.$$
\end{proof}

\begin{cor}\label{delta}
Let $\alpha$ be a locally constant $H$-labeling of G such that for all $h \in im(\alpha), \; val(h)=m$.
Define
\begin{enumerate}
\item $\Delta_1 = G \text{  and  } \alpha_1 = \alpha,$
\item $\pi_{n+1} : \Delta_{n+1} \rightarrow \Delta_{n}$ be the graph epimorphism induced by the first coordinate projection,
\item $\Delta_{n+1} = \Delta_n \zigzag_{\alpha_n} H \text{ and } \alpha_{n+1}=\alpha_n \circ \tilde{\pi}_{n+1},$
where $\tilde{ \pi }_{n+1}:D(\Delta_{n+1}) \rightarrow D(\Delta_n)$ is the map induced by ${\pi}_{n+1}$.
\end{enumerate}
Then $\pi_{n+1}$ is a combinatorial covering map for $n \geq 0.$ 
\end{cor}

Combining Corollary \ref{delta} with the result in Corollary \ref{laplace} we get:

\begin{cor}
For a finite graph $G$ and $H$-labeling satisfying the conditions of Corollary \ref{delta}, the spectrum of the normalized Laplacian for $\Delta_n$ is contained in the spectrum of the normalized Laplacian of $\Delta_m$ for all $m \geq n$.
\end{cor}

\end{document}